\documentclass[a4paper,10pt]{article}

\usepackage{latexsym}
\usepackage{amsfonts}
\usepackage[reqno]{amsmath}
\usepackage{amsthm}
\usepackage{amssymb}
\usepackage{bbm}
\usepackage{marvosym}
\usepackage{mathrsfs}
\usepackage{textcomp}

\usepackage[ps2pdf]{hyperref}

\newcounter{numtho}

\newtheorem{lemme}[numtho]{Lemma}
\newtheorem{corollaire}[numtho]{Corollary}
\newtheorem{proposition}[numtho]{Proposition}
\newtheorem{definition}[numtho]{Definition}
\theoremstyle{remark} 

\newcommand{\defe}[1]{\textbf{#1}}
\newcommand{\dpt}[3]{#1\,:\,#2\to #3}
\newcommand{\Dsddb}[4]{\frac{d}{d#2}\Big[#1\Big]_{#3=#4}}
\newcommand{\Dsdd}[3]{ \Dsddb{#1}{#2}{#2}{#3}   }
\def\emu{^{-1}}

\def\mfo{\vartheta}

\newcommand{\mapdp}[5]{%
                       \begin{equation}\begin{aligned}
                           #1\,:\,#2&\to #3\\
                                  #4&\to #5
                        \end{aligned}\end{equation}
                        }
\def\tq{\,|\,}
\def\fge{''}
\def\oge{``}

\def \ra{{\rightarrow}}
\def \lra{{\longrightarrow}}

\newcommand{\hS}{\mathscr{S}}
\newcommand{\hH}{\mathscr{H}}

\DeclareMathOperator{\pr}{\texttt{pr}}
\DeclareMathOperator{\id}{id}
\DeclareMathOperator{\tr}{Tr}
\DeclareMathOperator{\ad}{ad}
\DeclareMathOperator{\Ad}{Ad}
\DeclareMathOperator{\AD}{\textbf{Ad}}

\DeclareMathOperator{\End}{End}
\DeclareMathOperator{\Iso}{Iso}

\def\sA{\mathcal{A}}
\def\sG{\mathcal{G}}
\def\sH{\mathcal{H}}
\def\sK{\mathcal{K}}

\def\sN{\mathcal{N}}
\def\sP{\mathcal{P}}
\def\sQ{\mathcal{Q}}
\def\sR{\mathcal{R}}
\def\mG{\mathcal{G}}

\def\eR{\mathbb{R}}
\def\eZ{\mathbb{Z}}

\def\sldr{\mathfrak{sl}_2(\eR)}
\def\SLdr{SL_2(\eR)}
\def\soun{\mathfrak{so}(1,n)}

\def\sod{\mathfrak{so}(2)}
\def\son{\mathfrak{so}(n)}

   \def\dbar{\leavevmode\hbox to 0pt{\hskip.2ex
    \accent"16\hss}d}

\title{Black holes in symmetric spaces : anti-de Sitter spaces}
 \author{
{\bf L.~ Claessens} \thanks{``Chercheur FRIA'', Belgium}  \\{\it D\'epartement de Math\'ematique}\\ {\it
         Universit\'e Catholique de Louvain}\\ {\it
         Chemin du cyclotron, 2, }\\
         {\it 1348 Louvain-La-Neuve, Belgium}\\ E-mail:
        claessens@math.ucl.ac.be\\
        {\bf S.~Detournay} \thanks{``Chercheur FRIA'', Belgium}~
         \\{\it M\'ecanique et Gravitation}\\ {\it Universit\'e de Mons-Hainaut, 20
         Place du Parc}\\ {\it 7000 Mons, Belgium}\\E-mail:
        stephane.detournay@umh.ac.be}

\begin{document}
\maketitle

\begin{abstract}
Using symmetric space techniques, we show that closed orbits of
the Iwasawa
 subgroups of $SO(2,l-1)$ naturally define singularities of a black hole causal structure
  in anti-de Sitter spaces in $l \geq 3$ dimensions. In particular, we recover for $l=3$
   the non-rotating massive BTZ black hole. The method presented here is very simple and
    in principle generalizable to any semi-simple symmetric space.

\end{abstract}

\tableofcontents

\section{Introduction}

Causal black holes are distinguished from metric black holes by the fact
that they do not exhibit curvature singularities. They are
obtained as a quotient of certain spaces under the action of a
discrete isometry subgroup. To avoid closed time-like curves in
the resulting space, the parts of the original space where the
identifications are time-like must be cut out. In this context, the {\itshape raison d'\^etre} of the quotient operation
 is to make the resulting space \oge causally inextensible\fge.

The most celebrated examples, the BTZ black holes \cite{BTZ,BHTZ},
    are built from the three dimensional anti-de Sitter space ($AdS_3$) by identifying points
 along orbits of particular Killing vectors. They represent
 axisymmetric and static black hole solutions of (2+1)-dimensional
 gravity with negative cosmological constant. Furthermore, some of
 them enjoy remarkable Lie group-theoretical properties, pointed out in
 \cite{BDHRS,Keio}. The example which will be of interest for our
 purposes is the \emph{non-rotating massive} BTZ black hole. In
 this case, it was observed that the structure of the black hole
 singularities and horizons are closely linked to \emph{minimal
 parabolic} (or \emph{Iwasawa}) subgroups of $\widetilde{\SLdr} \simeq
 AdS_3$. We will mainly be concerned in this work to extend these
observations to higher-dimensional anti-de Sitter spaces.

Higher-dimensional generalizations of the BTZ construction have
been studied in the physics' literature, by classifying the
one-parameter isometry subgroups of Iso($AdS_l ) = SO(2,l-1)$, see
\cite{Figueroa,AdSBH,Madden,Banados:1997df,Aminneborg,HolstPeldan}.
Nevertheless, the approach we will adopt here is conceptually
different. We will first reinterpret the non-rotating BTZ black
hole solution using symmetric spaces techniques
and present an alternative way to express its singularities.
  The latter will be seen as closed orbits of Iwasawa subgroups of the isometry group.
As we will show, this construction extends straightforwardly to
higher dimensional cases, allowing to build a non trivial black
hole on anti-de Sitter spaces of arbitrary dimension $l\geq 3$.  A
groupal characterization of the event horizon is also obtained.
From this point of view, all anti-de Sitter spaces of dimension
$l\geq 3$ appear on an equal footing. For the sake of
completeness, we also analyze in some details in appendix B the
two-dimensional case, for which the construction does not yield a
black hole structure.

A natural question arising from this analysis is the following  : \emph{given
 a semi-simple symmetric space, when does the closed orbits of the Iwasawa subgroups
  of the isometry group seen as singularities define a non-trivial causal structure?}
   We here answer this question in the case of anti-de Sitter spaces, using techniques
    allowing in principle for generalization to any symmetric semi-simple symmetric space.

This paper is organized as follows. Section \ref{BTZ} is devoted
to the presentation of some aspects of the non-rotating BTZ black
holes. We state some properties whose proofs are left to appendix
A or to the existing literature. In the third section, we present
some general elements of the theory of symmetric spaces,
applicable to the study of anti-de Sitter spaces. We show how the
non-rotating BTZ black holes fit in this context and how the
singularities can be expressed in a way suitable for
generalization to higher dimensions. In section \ref{BHStructure},
we show that the proposed definition for the singularities indeed
gives rise to a black hole structure by proving the existence of
an event horizon, whose characterization is provided using the
Iwasawa decomposition of $\Iso(AdS_l)$. We leave the particular
two-dimensional case to appendix \ref{AdS2}, while some explicit
computation details are related in appendix \ref{app_calc}.

\section{BTZ black holes and minimal parabolic subgroups}\label{BTZ}


In this section, we recall for the reader's convenience the
definition and construction of the non-rotating BTZ black hole
\cite{BTZ,BHTZ}, emphasizing on some geometrical properties put
forwards in \cite{BRS,BDRS,Keio,MemClem}. To lighten the
presentation, the proofs will essentially be omitted and referred
to the existing literature, or recast in appendix \ref{AppBTZ}.

This situation will serve us as a guideline in defining black
holes in general anti-de Sitter spaces (see section
\ref{BHStructure}).

Ba\~nados, Henneaux, Teitelboim and Zanelli observed that taking
the quotient of (a part of) the three-dimensional anti-de Sitter
space ($AdS_3$) under the action of well-chosen discrete subgroups
of its isometry group gives rise to solutions which correspond to
axially symmetric and static black hole solutions of
(2+1)-dimensional Einstein gravity with negative cosmological
constant, characterized by their mass $M$ and angular momentum
$J$.

The space $AdS_3$ is defined as the (universal covering of the)
simple Lie group $\SLdr$

\begin{equation}\label{AdS3}
 AdS_3 \cong \SLdr = \{g \in GL_2(\eR) \, | \, \mbox{det} g = 1 \} := G
 \end{equation}
 endowed with its Killing metric $\dpt{B}{\sG\times\sG}{\eR}$,

\[
   B(X,Y)=\tr(\ad(X)\circ\ad(X))
\]
which can be extended to the whole group by

\begin{equation}
B_g(X,Y)=B(dL_{g\emu}X,dL_{g\emu}Y).
\end{equation}
Here, $\sG$ stands for the Lie algebra of $\SLdr$ :

\begin{equation}
\begin{split}
 \sG &:= \sldr = \{ X\in\End(2,\eR)\tq \tr(X) = 0\} \\
     &= \{ z^H H + z^E E + z^F F \}_{\{z^H,z^E,z^F \in \eR\}}.
\end{split}
\end{equation}

 The generators $H$, $E$ and $F$ satisfy the usual commutation
 relations

\begin{equation}
  [H,E] = 2E, \quad [H,F]=-2F, \quad [E,F]=H.
\end{equation}
We define the following one-parameter subgroups of $\SLdr$ :
\begin{equation}\label{IwasawaSubgroups}
A = \exp(\eR H), \quad  N = \exp(\eR E), \quad \bar{N}=\exp(\eR
F), \quad K=\exp(\eR T),
\end{equation}
with $T=E-F$. They are the building blocks of the Iwasawa
decomposition
\begin{equation}
 K \times A \times N \longrightarrow \SLdr : (k,a,n)
 \longrightarrow kan \quad \mbox{or} \quad ank.
 \end{equation}

The $\SLdr$ subgroups $AN$ and $A\bar{N}$ are called
\emph{Iwasawa subgroups}; they are minimal parabolic subgroups.

We will also use another representation of $AdS_3$, which can
equivalently be seen as the hyperboloid
\begin{equation} \label{hyperboloide}
 u^2 + t^2 - x^2 - y^2 = 1
 \end{equation}
embedded in $\eR^{2,2}$, that is the four-dimensional flat space
with metric $ds^2 = -du^2 - dt^2 + dx^2 + dy^2$.

>From \eqref{hyperboloide}, the isometry group of $AdS_3$, denoted by
$\Iso(G)$, is the four-dimensional Lorentz group $O(2,2)$. It
is locally isomorphic to $G \times G$, through the action
\begin{equation}
 (G \times G) \times G \longrightarrow G : ((g_L,g_R),z)
 \rightarrow g_l \, z g_R^{-1},
 \end{equation}
 which corresponds to the identity component of
 $Iso(G)$ (from the bi-invariance of the Killing
 metric), and because of the Lie algebra isomorphism
\begin{equation}\label{Iso}
 \Phi : \mG \times \mG \rightarrow iso(G) : (X,Y)
 \rightarrow  \overline{X} - \underline{Y},
\end{equation}
where $\overline{X}$ (resp. $\underline{Y}$) denotes the
right-invariant (resp. left-invariant) vector field on $G$
associated to the element $X$ (resp. $Y$) of its Lie algebra.

We have now dispose of all necessary ingredients to make the
definition of BTZ black holes more precise.

\begin{definition} \label{BHTZ}
The one-parameter subgroup of $Iso(G)$ defined by
\begin{equation}
 \psi_t (g) = \exp(t\, a\, H) \, g \exp(-t\, a\, H), \quad a \in \eR_0,\quad g \in G \,
 \end{equation}
 is called the \defe{BHTZ subgroup}. Its generator $ \Xi =
 a(H,H)$ is called the \defe{identification vector}.
 The \defe{BHTZ action} associated to $\Xi$ is $\psi_{\eZ} : G
 \rightarrow G$.
 \end{definition}

\begin{definition} \label{safe}
A \defe{safe region} in $AdS_3$ is defined as an open an connected
domain
\begin{equation}
 ||\Xi||^2 := \beta_z (\Xi,\Xi) > 0 \quad .
 \end{equation}
\end{definition}

\begin{definition}\label{NRMassive}
A \defe{non-rotating massive BTZ black hole} is obtained as the
quotient of a safe region in $AdS_3$ under the BHTZ action.
\end{definition}

This definition deserves some comments. First, the restriction to
a safe region in $AdS_3$ ensures that the resulting quotient space
be free of closed time-like curves. This means that the other
parts of $AdS_3$ have to be \oge cut out\fge{} from the original space.
Furthermore, due to the identifications, one may restrict to a
fundamental domain of the BHTZ action. Secondly, the \defe{black
hole singularities} $\hS$ are defined as the surfaces where
the identification vector becomes light-like :

\begin{equation}  \label{Singularities}
 \hS = \{z \in AdS_3 \tq \beta_z(\Xi,\Xi) = 0\}.
\end{equation}

 Thus, the BTZ black hole singularities represent singularities in
 the causal structure, not curvature ones. The resulting space is
 causally inextensible, i.e. trying to extend it would produce
 closed time-like curves. Finally, the BTZ space-time exhibits all
 characteristic features of a black hole. Namely, it has
 \emph{event horizons}, that is,
 surfaces hiding a
 region (\emph{the interior region}, see hereafter) causally
 disconnected from spatial infinity.

Note that it is the choice of identification vector which dictates
the nature (rotating, extremal, vacuum or non-rotating massive) of
the resulting black hole. Moreover, not all choices give rise to
black holes.

 The reason why we here focus on the non-rotating massive
case lies in the peculiar geometrical properties of its horizons
and singularities. To define the horizons properly, we will need
the concept of light-rays and light-cones issued from a point.

\begin{definition} \label{lightraycone}
A \defe{light-ray} starting from a point $g$ in a safe region
 is a
  curve

\begin{equation}
 l_g^k(s) = \exp(-s Ad(k) E)g
 \end{equation}
 for a given $k\in K$. The \defe{future} and \defe{past light-cones} at
 $g$ are given by

\begin{equation}
 C_g^{\pm} = \{l_g^k(s)\}_{%
\begin{subarray}{l}
k\in K \\ s\in \eR^{\pm}
 \end{subarray}}.
 \end{equation}
 \end{definition}

We are now ready to define the horizons.

\begin{definition}\label{interior-horizons}
A point $g$ will be said to lie in the \defe{future
interior region}, denoted by ${\cal M}^{\mbox{int},+}$, if all
future-directed light-rays issued from $g$ necessarily fall into
the black hole singularity, that is

\begin{equation}\label{interioreq}
 g \in {\cal M}^{\mbox{int},+} \Leftrightarrow \forall k \in K, \exists s \in
 \eR^+  \mbox{s.t.}  \|\underline{H} - \overline{H}\|^2_{l^k_g (s)}=0.
 \end{equation}
The \defe{future horizon} $\hH^+$ is defined as the boundary
of ${\cal M}^{\mbox{int},+}$.
\end{definition}
Equation \eqref{interioreq} simply expresses that any
future-directed causal signal necessarily falls into the black
hole singularity and cannot escape it. The \defe{past interior
region} and \defe{past horizon} are defined in a similar way.

Using the embedding \eqref{hyperboloide} of $AdS_3$ into $\eR^{2,2}$, one finds,
 from \eqref{Singularities} and \eqref{interioreq}, that
\begin{equation}\label{BTZSingHor}
\hS \equiv t^2 - y^2 = 0 \quad\textrm{and} \quad \hH \equiv u^2 - x^2 = 0,
 \end{equation}
 where $\hH = \hH^+ \cup \hH^-$.

These results can be stated more intrinsically as follows:
\begin{proposition}\label{BTZSing}
In $G=AdS_3$, the non-rotating BTZ black hole singularities are
given by a union of minimal parabolic subgroups of G:
\begin{equation}
\hS = Z(G) A N \cup Z(G) A \overline{N},
 \end{equation}
 where $Z(G)=\{e,-e\}$ denotes the center of $G=\SLdr$.
 \end{proposition}

\begin{proposition}\label{BTZHor}
In $G=AdS_3$, the non-rotating BTZ black hole horizons correspond
to a union of lateral classes of minimal parabolic subgroups of G
:
\begin{equation}
 \hH = Z(G) A N J \cup Z(G) A \overline{N} J \quad,
 \end{equation}
 where $J=\exp(\frac{3\pi}{2}T) \in K$ satisfies $J^2 = e$.
 \end{proposition}

These two propositions actually follows directly from
\eqref{BTZSingHor}, by using the parametrization $g=\left(
         \begin{array}{cc} u+x & y+t\\ y-t & u-x
         \end{array}
         \right)$.
They show that the black hole structure is closely related to the
minimal parabolic subgroup of $\SLdr$. Of course, this
construction cannot be generalized in a straightforward way to
higher-dimensional anti-de Sitter spaces, because of the peculiar
nature of the three-dimensional case, being the only to enjoy a
group manifold structure. Rather, we will reconsider in the next
section the case treated here in a more general framework, putting
on an equal footing all anti-de Sitter spaces. Again, a minimal
parabolic subgroup will reveal crucial in the construction.

\section{Symmetric space structure on anti-de Sitter}

Most of the material of this section can be found in a general
framework in \cite{Helgason, Loos, kobayashi, kobayashi2}.

\subsection{Basic facts}\label{BasicFacts}

As physical space, $AdS_l$ is the set of points
$(u,t,x_1,\ldots,x_{l-1})\in \eR^{2,l-1}$  such that
$u^2+t^2-x_1^2-\ldots-x_{l-1}^2=1$. The transitive (an isometric)
action of $SO(2,l-1)$ on $AdS_l$ yields an homogeneous space
structure. Let's parameterize the matrix representation of the
groups in such a way that $SO(1,l-1)$ --seen as a subgroup of
$SO(2,l-1)$-- leaves unchanged the vector $(1,0,\ldots,0)$. In
this case we have an homogeneous space isomorphism

\[
  AdS_l=\frac{SO(2,n)}{SO(1,n)}
\]
with $n=l-1$. The isomorphism is explicitly given by

\begin{equation}
[g]\to g\cdot
\begin{pmatrix}
1\\0\\\vdots
\end{pmatrix}
\end{equation}
where the dot denotes the action matrix times vector of the
 representant $g\in [g]$ in the defining representation of $SO(2,n)$.
  The classes are taken from the right :  $[g]=\{gh\tq h\in H\}$.

From now we set $G:=SO(2,n)$ and $H:=SO(1,n)$;
 the symbols $\sG$ and $\sH$ denote their respective
  Lie algebras.
    We also write $\mfo=[e]$ and $M=G/H=AdS_l$.
      We consider a Cartan involution $\dpt{\theta}{\sG}{\sG}$
       which gives a Cartan decomposition

\[
\sG=\sK\oplus\sP,
\]
and an involutive automorphism $\sigma=\id|_{\sH}\oplus(-\id)|_{\sQ}$
 which gives a reductive symmetric space decomposition

\[
\sG=\sH\oplus\sQ
\]
with

\begin{equation}
[\sH,\sH]\subset\sH,\quad[\sH,\sQ]\subset\sQ,\quad [\sQ,\sQ]\subset\sH.
\end{equation}
One can choose them in such a manner that $[\sigma,\theta]=0$.

 The space $\sQ$ can be identified with the tangent space
 $T_{[e]}M$. We can extend this identification by defining
 $\sQ_g=dL_g\sQ$. In this case $\dpt{d\pi}{\sQ_g}{T_{[g]}M}$
  is a vector space isomorphism.

 The last point is to find Iwasawa decompositions
  $\sH=\sA_{\sH}\oplus\sN_{\sH}\oplus\sK_{\sH}$ and
   $\sG=\sA\oplus\sN\oplus\sK$ with $\sA_{\sH}\subset\sA$ and
    $\sN_{\sH}\subset\sN$. We denote by $A$, $N$ and $K$ the
     exponentials of $\sA$, $\sN$ and $\sK$; and $\overline{N}=\theta(N)$.

Some explicit matrix
   choices are given in appendix \ref{app_calc}.
Since the Killing form $B$ is an $\Ad_H$-invariant product on
$\sQ$, we can define

\begin{equation}
B_g(X,Y)=B_g(dL_{g\emu}X,dL_{g\emu}Y)
\end{equation}
which descent to an homogeneous metric on $T_{[g]}M$ :

\begin{equation}
B_{[g]}(d\pi X,d\pi Y)=B_g(\pr X,\pr Y)
\end{equation}
where $\dpt{\pr}{T_gG}{dL_g\sQ}$ is the canonical projection. Properties of this product are given in \cite{Kerin}.

\subsection{Causal structure on anti-de Sitter space}      \label{Causal}

Let us start this section by computing the closed orbits of the
action of $AN$ and $A\overline{N}$ on $AdS_l$. In order to see if
$x=[g]\in M$ lies in a closed orbit of $AN$, we \oge compare\fge{}
the basis $\{d\pi dL_g q_i\}$ of $T_xM$ and the space spanned by
the fundamental vectors of the action. If these two spaces are the
same then $x$ belongs to an open orbit (because a submanifold is
open if and only if it has same dimension as the main manifold). This
idea is precisely contained in the following proposition.

\begin{proposition}   \label{prop:pr_ouvert}
If $R$ is a subgroup of $G$ with Lie algebra $\sR$, then the orbit $R\cdot \mfo$ is open in $G/H$ if and only if the projection $\dpt{\pr}{\sR}{\sQ}$ is surjective.
\end{proposition}

In order to check the openness of the $R$-orbit of $[g]$,
 we look at the openness of the $\AD(g\emu)R$-orbit of $\mfo$ using the proposition.

A great simplification is possible.
The $AN$-orbits are trivially $AN$-invariant.
 So the $K$ part of $[g]=ank$ alone fix the orbit in which $[g]$ belongs.
  In the explicit parametrization of $K$, we know that the $SO(n)$ part is
  \oge killed\fge{} by the quotient with respect to $SO(1,n)$.
  In definitive, we are left with \emph{at most} one $AN$-orbit for each element in $SO(2)$.
   Computations using proposition \ref{prop:pr_ouvert} show that the closed orbits are given
    by

\begin{equation}\label{Sing2}
\hS=\{\pm[AN],\pm[A \overline{N}]\}.
\end{equation}

We are now in position to make a link with the non-rotating BTZ
black hole we discussed in the previous section, through the
following

\begin{proposition}\label{BTZOpenOrbits}
The singularities of the non-rotating BTZ black hole, given in
\eqref{BTZSingHor}, coincide with the closed orbits of the action of
the subgroups $AN$ and $A \overline{N}$ of $SO(2,2)$ on $AdS_3$.
\end{proposition}
This may be checked by computing the fundamental vector fields of
the actions of $AN$ and $A \overline{N}$, then by determining the
loci where they span a space of dimension less than 3, and finally
observing that this actually corresponds to the equation for $\hS$
in \eqref{BTZSingHor}.


The advantage of this reinterpretation is that it allows, this
time, for a straightforward generalization to higher-dimensional
anti-de Sitter spaces. Proposition \ref{BTZOpenOrbits} motivates the
following

\begin{definition}\label{Singular}
A point in $AdS_l$ is \defe{singular} if it belongs to a closed
orbit of the Iwasawa group $AN$ or $A\overline{N}$.
\end{definition}

This definition finds its origin in the next proposition, which we
will mainly concerned with in the next section :

\begin{proposition}\label{LeBut}
In $AdS_l$, for $l \geq 3$, defining singularities as the closed
orbits of the Iwasawa subgroups $AN$ and $A\overline{N}$ of
$SO(2,l-1)$ gives rise to a black-hole structure, in the sense
that there exists a non empty event horizon.
\end{proposition}

Let us make this more precise. As in the three dimensional case,
we need to define the notion of light-cone in $AdS_l$.

General theory about symmetric spaces says that if $E$ is
nilpotent in $\sQ$, then  $\{\Ad(k)E\}_{k\in K_H}$ is the set of
all the light-like vectors in $T_{[\mfo]}AdS_l\simeq\sQ$. So the
future light cone of $\mfo$ is given by

\[
  C^+_{[\mfo]}=\{\pi(e^{-t\Ad(k)E})\}_{%
\begin{subarray}{l}
t\in\eR^+\\k\in K_H
\end{subarray}}
\]
and the one of a general element $[g] \in AdS_l$ is obtained by the (isometric)
 action of $g$ thereon :

\begin{equation}  \label{eq:exprcone}
C^+_{\pi(g)}=\{\pi(ge^{-t\Ad(k)E})\}_{%
\begin{subarray}{l}
t\in\eR^+\\k\in K_H
\end{subarray}}.
\end{equation}
It should be noted that this definition is independent of the
choice of the representant $g$ in the class $\pi(g)$ because, for
any $h\in H$, $\pi(ghe^{-t\Ad(k)E})=\pi(ghe^{-t\Ad(k)E}h\emu)$
which is simply a reparametrization in $K_H$.

We are now able to define the causality as follows.
 A point $[g]\in AdS_l$ belongs to the \defe{interior region} if for all direction $k\in K_H$,
  the future light ray through $[g]$ intersects the singularity within a \emph{finite} time.
   In other words, it is interior when the whole light cone ends up in the singularity.
    A point is \defe{exterior} when it is not interior.
 A particularly important set of point is the \defe{event horizon}, or simply \emph{horizon}, defined as the boundary
  of the interior. When a space contains a non trivial causal structure
   (i.e. when there exists a non empty horizon),
    we say that the definition of singularities gives rise to a \defe{black hole}.

\section{Black hole structure on anti-de Sitter spaces} \label{BHStructure}

\subsection{General method for computing the singularities}

First, let us give an alternative to proposition
\ref{prop:pr_ouvert}
 to study the openness of an $AN$-orbit. We denote by $\hS_{AN}$ the closed orbits
  of $AN$ and by $\hS_{A\overline{N}}$ the ones of $A\overline{N}$.
   We explain the method for $\hS_{AN}$, but the same with trivial adaptations
    is true for $\hS_{A\overline{N}}$.

If $x\in M$ belongs to $\hS_{AN}$, the tangent space of his $AN$-orbit has lower
 dimension that the tangent space of $M$.
  In this case the volume spanned by the fundamental vectors at $x$ is zero.
 The idea is simple : we build the volume form
 $\nu_x$ of $T_xM$ and we apply it on a basis of the fundamental fields.
 If the result is zero, then $x$ belongs to the $\hS_{AN}$.
 The action is given by

\mapdp{\tau}{AN\times M}{M}{(an,[g])}{[ang].}
If $X\in\sA\oplus\sN$ and $[g]\in M$, then

\begin{equation}
  X^*_{[g]}=-d(\pi\circ R_g)X.
\end{equation}
As said before, if $\{q_i\}$ is a basis
 of $\sQ$ then a basis of $T_{[g]}M$ is given by $\{d\pi dL_gq_i\}$. We define

\[
\nu=q_0^{\flat}\wedge q_1^{\flat}\wedge \ldots \wedge
q_{l-1}^{\flat}
\]
where $q_{i[g]}^{\flat}=B_{[g]}(d\pi dL_g q_i,\cdot)$. The
condition for $[g]$ to belongs to $\hS_{AN}$ reads

\begin{equation}\label{eq:nusurN}
\nu_{[g]}(N_1^*{}_{[g]},N_2^*{}_{[g]},\ldots,N_l^*{}_{[g]})=0
\end{equation}
for all choices of $N_j$ in a basis of $\sA\oplus\sN$. It
corresponds to the vanishing of $l \times l$ determinants. Our
purpose is now to compute the products

\[
\begin{split}
  B_{[g]}(d\pi dL_g q_i,N^*_j{}_{[h]})&=-B_g(dL_g q_i,dR_g N_j)\\
                   &=-B_e(q_i,\Ad(g\emu)N_j).
\end{split}
\]
We note
\[
\Delta_{ij}([g])=B(q_i,\Ad(g\emu)N_j)
\]
where $N_j$ runs over a basis of $\sA\oplus\sN$ and $q_i$ a one of
$\sQ$. Our problem of light cone (see explanations around
expression \eqref{eq:exprcone}) leads us to compute

\begin{equation} \label{eq:elemtr}
\Delta_{ij}(\pi(ge^{-tk\cdot E}))=B(\Ad(e^{-tk\cdot
E})q_i,\Ad(g\emu)N_j)
\end{equation}
where $k\cdot E$ is a notation for $\Ad(k)E$.

A way to proceed is to express all our elements of $SO(2,n)$ in the
root space decomposition
\[
\sG=\sG_{(0,0)}\bigoplus_{\lambda\in\Sigma}\sG_{\lambda}.
\]
The purpose
 of this resides in the fact
  that the Killing form $B(X,Y)$ is most easy to compute when $X$ and $Y$ are
   in some root spaces. In order to be more synthetic in the text, most of explicit
  decompositions are given in appendix \ref{app_calc}.

An important computational remark is the fact that $E$ is
nilpotent, so $\Ad(k)E$ also is and
$\Ad(e^{-t\Ad(k)E})X=e^{-t\ad(k)E}X$ only gives second order
expressions with respect to $t$. These computations are
nevertheless heavy, but can fortunately be circumvented by a
simple counting of dimensions, as we describe in the next
subsection.

\subsection{$AdS_l$-adapted method for computing the singularities}

We here explicitly use the description of $AdS_l$ in terms of the
embedding coordinates $(u,t,x_1,\ldots,x_{l-1})\in \eR^{2,l-1}$
(see \autoref{BasicFacts}), and the choices of generators related
in appendix \ref{app_calc} .

\begin{proposition}
In term of the embedding of $AdS_l$ in $\eR^{2,l-1}$, the closed
orbits of $AN \subset SO(2,l-1)$ are located at $y-t = 0$.
Similarly, the closed orbits of $A \overline{N}$ correspond to
$y+t=0$.
\end{proposition}

In other words, the equation
\begin{equation}
t^2-y^2=0
\end{equation}
describes the singularity $\hS=\hS_{AN}\cup\hS_{A\overline{N}}$.

\begin{proof}
The different fundamental vector fields of the $AN$ actions can be
computed by $X^*_{[g]}=-Xg\cdot\mfo$. For example, in $AdS_3$,

\[
\begin{split}
   M^*_{[g]}&=
\begin{pmatrix}
0&-1&0&1\\
1&0&-1&0\\
0&-1&0&1\\
1&0&-1&0
\end{pmatrix}
\begin{pmatrix}
u\\t\\x\\y
\end{pmatrix}
=
\begin{pmatrix}
-t+y\\u-x\\-t+y\\u-x
\end{pmatrix}\\
&=(y-t)\partial_u+(u-x)\partial_t+(y-t)\partial_x+(u-x)\partial_y.
\end{split}
\]
Full results are

\begin{subequations}\label{Gen}
\begin{align}
J_1^*&=-y\partial_t-t\partial_y\\
J_2^*&=-x\partial_u-u\partial_x                                                      \label{eq:Jds}\\
M^*  &=(y-t)\partial_u+(u-x)\partial_t+(y-t)\partial_x+(u-x)\partial_y\\
L^*  &=(y-t)\partial_u+(u+x)\partial_t+(t-y)\partial_x+(u+x)\partial_y\\
W_i^*&=-x_i\partial_t-x_i\partial_y+(y-t)\partial_i\\
V_j^*&=-x_j\partial_u-x_j\partial_x+(x-u)\partial_j \quad ,\quad
i,j=3,\ldots,l-1 \label{eq:Vjs}
\end{align}
\end{subequations}

First consider points satisfying $t-y=0$. It is clear that, at
these points, the $l$ vectors $J_1^*$, $M^*$, $L^*$ and $W_i^*$
are linearly dependent. Then, there are at most $l-1$ linearly
independent vectors amongst the $2(l-1)$ vectors \eqref{Gen}, thus
the points belong to a closed orbit.

We now show that a point with $t-y\neq 0$ belongs to an open orbit
of $AN$. It is easy to see that $J_1^*$, $L^*$ and $M^*$ are three
linearly independent vectors. The vectors $V_i^*$ gives us $l-3$
more. Then they span a $l$-dimensional space.

The same can be done with the closed orbits of $A\overline{N}$.
The result is that a points belongs to a closed orbit of
$A\overline{N}$ if and only if $t+y=0$.
\end{proof}

\begin{corollaire}
The singularities coincide with the set of points in $AdS_l$ where
$\| J_1^* \|^2 = 0$.
\end{corollaire}
This generalizes proposition \ref{BTZOpenOrbits} to any dimension. Hence, a discrete quotient of $AdS_l$ along orbits of $J_1^*$ gives a direct higher-dimensional generalization of the non-rotating BTZ black hole.

\subsection{Existence of an horizon}
Let us show that the definition \ref{Singular} gives rise to a
black hole causal structure, namely that it leads to the existence
of horizons, as defined in \autoref{Causal}.

We first consider points of the form $K\cdot\mfo$, which are
parameterized by an angle $\mu$.
Up to choice
of this parametrization,
 a light-like geodesic trough $\mu$ is given by
 \begin{equation}
  K\cdot \mbox{e}^{-s\Ad(k)E}\cdot\mfo
\end{equation}
with $k\in SO(l-1)$ and  $s\in\eR$.

This geodesic reaches $\hS_{AN}$ and $\hS_{A\overline{N}}$ for
values $s_{AN}$ and $s_{A\overline{N}}$ of the affine parameter,
given by

\begin{equation}   \label{eq:tempssingul}
 s_{AN} = \frac{\sin\mu}{\cos\mu - \cos\alpha},\quad \text{and} \quad  s_{A\overline{N}} = \frac{\sin\mu}{\cos\mu + \cos\alpha}
\end{equation}
where $\cos\alpha$ is the second component of the first column of $k$, see appendix \ref{app_calc} and equation \eqref{eq:AdkE}.

Because the part $\sin \mu =0$ is $\hS_{AN}$, we may restrict
ourselves to the open connected domain of $AdS_l$ given by $\sin
\mu > 0$. More precisely, $\sin\mu=0$ is the equation of
$\hS_{AN}$ is the $ANK$ decomposition. In the same way,
$\hS_{A\overline{N}}$ is given by $\sin\mu'=0$ in the
$A\overline{N}K$ decomposition.
 In order to escape the singularity, the point $\mu$ needs $s_{AN},s_{A\overline{N}}<0$.
  It is only possible to find directions (i.e. an angle $\alpha$)
   which gives it when $\cos u<0$.
   So the point $\cos u=0$ is one point of the horizon.

This proves proposition \ref{LeBut}. Remark that the two-dimensional case here appears
 as degenerate. Therefore, it is treated in appendix \ref{AdS2}, where we show that
  \emph{no black hole arises from this construction in $AdS_2$}.

\subsection{A characterization of the horizon}

Let $D[g]$ be the set of the light-like directions (vectors in $SO(n)$)
 for which the point $[g]$ falls into $\hS_{AN}$.
Similarly, the set $\overline{D}[g]$ is the one of directions
which fall into $\hS_{A\overline{N}}$. A great result is the fact
that it is possible to express $\overline{D}$ in terms of $D$.
Indeed

\begin{equation}
\begin{split}
k\in \overline{D}[g]&\textrm{ iff }\pi(ge^{tk\cdot E})\in\hS_{A \overline{N}}\\
                    &\textrm{ iff }\pi(\theta(g)\theta(e^{tk\cdot E_1}))\in\hS_{AN}\\
                    &\textrm{ iff }\theta(k)\in D(\theta[g])\\
                    &\textrm{ iff }k\in(D(\theta[g]))_{\theta}.
\end{split}
\end{equation}
So
\begin{equation} \label{eq:DbarD}
\overline{D}[g]=(D\theta[g])_{\theta}
\end{equation}
 where the definition of $k_{\theta}$ is

\[
\theta(\Ad(k)E)=\Ad(k_{\theta})E.
\]
This definition is possible because $\theta$ is an inner automorphism.

 It is easy to see that $\theta$ changes the sign of the spatial
  part of $k$, i.e. changes $w_i\to -w_i$.

How to express the condition $g\in\hH$ in terms of $D[g]$ ? The
condition to be in the black hole is $D[g]\cup
\overline{D}[g]=SO(n)$. If the complementary of $D[g]\cup
\overline{D}[g]$ has an interior (i.e. if it contains an open
subset), then by continuity the complementary $D[g']\cup
\overline{D}[g']$ has also an interior for all $[g']$ near $[g]$.
In this case, $[g]$ cannot belong to the horizon. So a
characterization of $\hH$ is the fact that the boundary of $D[g]$
and $\overline{D}[g]$ coincide. Equation \eqref{eq:DbarD} shows
that $\hH$ is $\theta$-invariant.

We can explicitly express $D[u]$ for $u\in SO(2)$ by examining
equation \eqref{eq:tempssingul}. Let us write $w_2$ instead of
$\cos \alpha$. The set $D[u]$ is the set of $w_2\in [-1,1]$ such
that $\cos u-w_2>0$ :

\begin{equation}
  D[u]=[-1,\cos\mu[ .
\end{equation}
So in order for $[u]$ to belong to $\hH$, it must satisfy

\[
  D[\theta]_{\theta}=[-1,\cos\mu'[_{\theta}=]-\cos\mu',1].
\]
Consequently, if $u$ is the $K$ component of $g$ in the $ANK$ decomposition and $u'$ the one of $\theta u$, then we can describe the horizon by

\begin{equation}
\cos u=-\cos u'
\end{equation}

\section*{Acknowledgments}
We are grateful to Pierre Bieliavsky for suggesting the problem. We also thank him, as well as Philippe Spindel, for numerous enlightening discussions. We would also like to thank
the ``Service de Physique Th\'eorique'' of the Universit\'e Libre de
Bruxelles, where part of this work has been achieved, for its
hospitality.

\appendix

\section{Horizons of the non-rotating BTZ black holes}\label{AppBTZ}

\subsection{Global description of the black hole}

In this appendix, we use results and techniques of
\cite{BRS,BDHRS,Keio,MemClem} to derive the equation of the
non-rotating BTZ black holes horizons. We will begin by stating
some results which will be useful in describing the global
geometry of the black hole.

\begin{proposition} Let $\sigma$ be the unique exterior
automorphism of $G$ fixing pointwise the Cartan subgroup $A$ and
consider the following \emph{twisted} action of $G$ on itself :
\begin{equation}
 \tau : G \times G \lra G : (g,x) \ra \tau_g(x) :=
 g\,x\,\sigma(g^{-1}).
 \end{equation}
Then, the BHTZ action (see definition \ref{BHTZ}) can be rewritten as
\begin{equation}
 \psi_n = \tau_{\exp(n \sqrt{M} H)}, \quad n\in \eZ.
 \end{equation}
\end{proposition}
The proof follows from the fact that $\sigma$ fixes the generator
$H$. Using the action $\tau$, one finds the following global
decomposition of $G$ :
\begin{proposition}
 The map
 \begin{equation}
 \phi : A \times G/A \lra G : (a,[g]) \ra \phi(a,[g]) := \tau_g(a)
 \end{equation}
 is well-defined as a global diffeomorphism.
 \end{proposition}

 This follows from the observation that the application
 \begin{equation} \label{twistedI}
 \phi : K \times A \times N \rightarrow G : (k,a,n) \rightarrow
 \phi(k,a,n) = \tau_{kn}(a) \quad
 \end{equation}
  is a global diffeomorphism on $G$ (``twisted Iwasawa
  decomposition''). As a consequence, the space $G$ appears as the
  total space of a trivial fibration over $A = SO(1,1) \simeq \eR$
  whose fibers are the $\tau_G-$orbits, i.e. the $\sigma-$twisted
  conjugacy classes. As a homogeneous $G-$space, every fiber is isomorphic to
  $G/A = AdS_2$. Moreover, the BHTZ action is fiberwise,
  because
  \begin{equation} \label{fiberwise}
 \tau_h (\phi(a,[g])) = \phi(a,h.[g]) = \phi(a,[hg]).
 \end{equation}
The Killing metric on $G$ turns out to be globally diagonal with
respect to the twisted Iwasawa decomposition \cite{BRS} :
\begin{equation}
ds^2_G = da_A^2 - \frac{1}{4} \cosh^2 (a) ds^2_{G/A},
\end{equation}
where $ds^2_{G/A}$ denotes the canonical projected $AdS_2-$metric
on $G/A$. The study of the quotient space $G/\eZ$ therefore
reduces to the study of $(G/A)/\eZ$.

The space $G/A$ can be realized as the $G-$equivariant universal
covering space of the adjoint orbit ${\cal O}:= \mbox{Ad}(G)H$ in
$\sldr$, where it corresponds to a one sheet hyperboloid. In this
picture, we may identify the part of the hyperboloid corresponding
to a safe region (see definition \ref{safe}) in $G$.

\begin{lemme}
In ${\cal O}$, a connected region where the orbits of the BHTZ
action are space-like is given by
\begin{equation}
 \{ X = x^H H+x^E E + x^F F \in {\cal O} \,\, | \,\, -1 < x^H <
 1\}.
\end{equation}
Furthermore, it can be parameterized as
\begin{equation}
 X = Ad\left(\exp(\frac{\theta}{2}H) \, \exp(-\frac{\tau}{2}
 T)\right)\, H, \, 0<\tau<\pi \, ,\,  -\infty<\theta<+\infty
 .
\end{equation}
\end{lemme}
This has been proven in \cite{BRS}. From this and the preceding
proposition, we find a global description of a safe region in $G$
well adapted to the BHTZ identifications, through the

\begin{proposition}
A global description of a safe region in $G$ is given by
\begin{equation}\label{CoordGlob}
 z(\rho,\theta,\tau) =
\tau_{\exp(\frac{\theta}{2}H) \,
 \exp(-\frac{\tau}{2}
 T)}(\exp(\rho H)).
 \end{equation}
 Furthermore, the action of the BHTZ subgroup reads in these coordinates
\begin{equation}
 (\tau , \rho, \theta) \rightarrow (\tau, \rho, \theta + 2 n
a).
\end{equation}
\end{proposition}

\subsection{Derivation of the horizons}

We now have to study the equation of \eqref{interioreq}. Using the
bi-invariance of the Killing metric and the Ad-invariance of the
Killing form, it reduces to
\begin{equation}\label{EqHoriz}
B(H,H) - B(H,Ad(\mbox{e}^{-s Ad(k)E}) Ad(x) H) = 0.
\end{equation}

\begin{lemme}\label{BiInv}
${\cal M}^{\mbox{int},+}$ is A bi-invariant.
\end{lemme}
 \begin{proof}
 This equation is clearly invariant under $x \to x\cdot a \, , \, a\in
A$. In order to see the invariance under $x \to a\cdot x$, one uses the cyclicity
 of the trace to
bring the second term to

\[
 B(H,Ad(Ad(a^{-1})\mbox{e}^{u Ad(k)E}) Ad(x) H).
\]
But $Ad(a^{-1})\mbox{e}^{-s Ad(k)E} =
\mbox{e}^{-\tilde{s} Ad(\tilde{k})E}$, with $\tilde{s} = s
(\mbox{e}^{-2a} \cos^2\theta + \mbox{e}^{2a} \sin^2 \theta)$ and
$\cot t = \mbox{e}^{-2a} \cot \theta$, where $k = \mbox{e}^{\theta
T}$ and $\tilde{k}=\mbox{e}^{t T}$. The net result is thus simply
a relabelling of the parameters (note that $s$ and $\tilde{s}$
have the same signs!)
\end{proof}

Let us now consider a light-ray (definition \ref{lightraycone}) starting from a
safe region in $G$. Because of \eqref{CoordGlob} and lemma
\ref{BiInv}, we may restrict our study to
\begin{eqnarray}
 z &=& \mbox{e}^{-\tau/2 T} \mbox{e}^{\rho H} \sigma(\mbox{e}^{\tau/2 T}) \\
 &=& \mbox{e}^{-\tau/2 T} \mbox{e}^{\rho H} \mbox{e}^{-\tau/2 T}.
 \end{eqnarray}
The equation to study reduces to
\begin{equation}\label{EqHoriz2}
B(H,H) - B(H,Ad(\mbox{e}^{-s Ad(k)E}) Ad(\mbox{e}^{-\tau/2 T}
\mbox{e}^{\rho H} \mbox{e}^{-\tau/2 T}) H) = 0
\end{equation}
with $\tau\in \, ]0,\pi[$ and  $\rho \in \eR$.

Let us focus on the points in $Ad(G)H$ corresponding to

\[
{\cal B} := Ad(\mbox{e}^{-\tau/2 T} \mbox{e}^{\rho H} \mbox{e}^{-\tau/2 T})H,
\]
with $\tau \in \, ]0,\pi[ \,\, , \,\, \rho \in \eR$. First note
that $Ad(\mbox{e}^{\rho H} \mbox{e}^{-\tau/2 H}) H$ precisely
corresponds to a safe region on the hyperboloid . Thus ${\cal B}$
is the region swept out by the a safe region when rotating it
counterclockwise around the $T$-axis with an angle~$\pi$.

It can be seen that the domain ${\cal B}$ can be decomposed into
three regions :

\begin{equation}
 {\cal B} = {\cal B}_1 \cup {\cal B}_2 \cup {\cal B}_3,
 \end{equation}
with

\begin{subequations}
\begin{align}
{\cal B}_1 &=Ad(A)Ad(\mbox{e}^{-\beta/2 T})H &&   \beta \in \, ]0,2\pi[,\\
{\cal B}_2 &=Ad(A)Ad(\mbox{e}^{t (E+F)})H    &&    t\in \eR, \\
{\cal B}_3 &=Ad(A)(-H \pm E) \textrm{ or } Ad(A)(-H \pm F).
\end{align}
\end{subequations}
Thanks to the A bi-invariance, we may forget about the $Ad(A)$ in
the above equations. We are thus led to analyze
the existence of solutions of \eqref{EqHoriz2} with $X \in {\cal B}$
of the form $X_1 = Ad(\mbox{e}^{-\beta/2 T})H$, $X_2
=Ad(\mbox{e}^{t (E+F)})H$ and $X_3 =-H \pm E \, , \, -H \pm F$.

Consider the first case. With $ Ad(\mbox{e}^{-\tau/2 T}
\mbox{e}^{\rho H} \mbox{e}^{-\tau/2 T}) H$ of the form
$Ad(\mbox{e}^{-\beta/2 T})H$, \eqref{EqHoriz2} becomes the following
equation :

\begin{equation}  \label{EqU}
 \frac{1}{4} s^2 (\cos \beta - \cos(\beta + 4\theta)) + s \sin
 \beta  + 2 \sin^2 \beta = 0.
\end{equation}
We are looking for the values of $\beta$ for which this equation
admits a solution for $s>0$, for all $\theta \in \, [0,\pi]$ --this range for $\theta$ originates from the fact that $G/A$ is a $\eZ_2$ covering of $\Ad(G)H$.
By considering the particular case $\theta = 0$, we find $s=-\tan
\frac{\beta}{2}$, thus the allowed values of $\beta$ have to lie
in the range $]\pi,2\pi[$. Let us look at the constrains imposed by
other values of $\theta$. We denote by $s_1$ and $s_2$ the two
roots of \eqref{EqU}. We have

\begin{eqnarray}
 s_1\cdot s_2 &=& \frac{4 \sin^2 \beta/2}{\sin 2\theta \sin(\beta +
 2 \theta)}, \label{PrRac} \\
 s_1 + s_2 &=& \frac{- 2 \sin \beta}{\sin 2\theta \sin(\beta +
 2\theta)}.
 \end{eqnarray}
First note that, $\forall \beta \in ]0,2\pi[$, $\sin 2\theta
\sin(\beta + 2\theta)$ may be positive or negative as $\theta$
varies in the range $[0,\pi]$. If $\sin \beta < 0$, then there are
two positive roots when $\sin 2\theta \sin(\beta + 2\theta)>0$,
and one positive and one negative when $\sin 2\theta \sin(\beta +
2\theta)<0$. Thus there always exist a positive solution for $u$,
for any $\theta$. If $\sin \beta > 0$, there are two negative
roots when $\sin 2\theta \sin(\beta + 2\theta)>0$. Consequently,
the interior region will correspond to points $X_1 =
Ad(\mbox{e}^{-\beta/2 T})H \, \, , \, \, \beta \in \,]\pi , 2\pi[
$ on the adjoint orbit.

For the second case, $X_2 =Ad(\mbox{e}^{t (E+F)})H$, the equation
we get is
\begin{equation}
 \frac{1}{4} s^2 (\cosh 2t - \cos 4\theta \cosh 2t + 2 \sin
 2\theta \sinh 2t) + s \cos 2\theta \sinh 2t + (1 - \cosh 2t) = 0.
 \end{equation}
 By considering two special cases, it is easy to see that this
 equation does not admit a positive solution in $u$ for all
 $\theta$. Indeed, for $\theta = \pi/2$, one finds $s =-\tanh t$,
 while for $\theta = 0$, one gets $s=\tanh t$. Thus there is no
 $t\ne 0$ satisfying both conditions. The last case yields no
 positive solution for all $\theta$ neither.

 As a conclusion we find that the interior region is given by
 \begin{equation}
  x \in {\cal M}^{\mbox{int},+} \Leftrightarrow Ad(x)H =
  Ad(A)Ad(\mbox{e}^{-\beta/2 T})H, \quad \mbox{with} \,\,
  \beta \in \, ]\pi , 2\pi[.
  \end{equation}
  The boundaries of the corresponding region in $Ad(G)H$ are given by $-H +
  r^2 E$ and $-H + r^2 F$ or
  \begin{equation}
   Ad(N^-)(-H) \cup  Ad(\bar{N}^+)(-H),
   \end{equation}
   with $N^- = \{\mbox{e}^{tE}\}_{t\leq 0}$ and $\bar{N}^+ = \{\mbox{e}^{tF}\}_{t\geq
   0}$.

The horizons can be deduced as
\begin{equation}
 x \in {\cal H^+} \Leftrightarrow Ad(x)H = Ad(N^-)(-H) \,\,
 \mbox{or} \,\, Ad(x)H = Ad(\bar{N}^+)(-H).
 \end{equation}
Because of the A-invariance, we may write
$x=\tau_{\mbox{e}^{-\frac{\tau}{2} T}}(\mbox{e}^{\rho H})$ and
look for the relation between $\tau$ and $\rho$ such that
\begin{equation}
 Ad(\tau_{\mbox{e}^{-\tau/2
T}}(\mbox{e}^{\rho H}))H = Ad(N^-)Ad(\mbox{e}^{\pi/2 T}) H.
\end{equation}
This amounts to require that
\begin{equation}\label{CondHor}
 \left(\mbox{e}^{-\tau/2 T} \mbox{e}^{\rho H} \mbox{e}^{-\tau/2 T}\right)^{-1} \,
 (\mbox{e}^{-t^2 E} \mbox{e}^{\pi/2 T}) \in A \cup Z(G).
 \end{equation}

This condition gives $\cos\tau = \tanh \rho$, $\rho<0$, $\tau\in
\, ]\pi/2,\pi[$. By replacing $\mbox{e}^{-t^2 E}$ with
$\mbox{e}^{t^2 F}$, one gets $\cos\tau = -\tanh \rho$, $\rho>0$,
$\tau\in \, ]\pi/2,\pi[$.

The domain ${\cal M}^{\mbox{int},-}$ is of course defined as
\begin{equation}\label{interior}
 x \in {\cal M}^{\mbox{int},-} \Leftrightarrow \forall k \in K, \exists u \in
 \eR^-  \mbox{s.t.}  \|\underline{H} - \overline{H}\|^2_{l^k_x (u)}=0.
 \end{equation}
The \emph{past horizon} ${\cal H^-}$ is defined as the boundary of
${\cal M}^{\mbox{int},-}$. By proceeding the same way, we find
that
\begin{equation}
  x \in {\cal M}^{\mbox{int},-} \Leftrightarrow Ad(x)H =
  Ad(A)Ad(\mbox{e}^{-\beta/2 T})H, \quad \mbox{with} \,\,
  \beta \in \, ]\pi,2 \pi[,
  \end{equation}
and
\begin{equation}
 x \in {\cal H^-} \Leftrightarrow Ad(x)H = Ad(N^+)(-H) \,\,
 \mbox{or} \,\, Ad(x)H = Ad(\bar{N}^-)(-H),
 \end{equation}
or in coordinates : $\tau\ \in \, ]0,\pi/2[$, $\cos \tau =
\tanh \rho$ for $\rho > 0$ and $\cos \tau =-\tanh \rho$ for $\rho <
0$.

We thus established the following
\begin{proposition}
In a safe region in $G$ parameterized by

\[
z(\rho,\theta,\tau) =\tau_{\exp(\frac{\theta}{2}H) \, \exp(-\frac{\tau}{2} T)}(\exp(\rho H)),
\]
 the horizons ${\cal H} := {\cal H}^+ \cup {\cal H}^-$ of the non-rotating BTZ black hole are given by
\begin{equation}
 \cos \tau = \pm \tanh \rho.
 \end{equation}
\end{proposition}
As a direct consequence, we have the
\begin{corollaire}
In terms of the embedding coordinates \eqref{hyperboloide} of $G$ in
$\eR^{2,2}$, the horizons of the non-rotating BTZ black hole are
\begin{equation}
 {\cal H} \equiv u^2 - x^2 = 0.
 \end{equation}
 \end{corollaire}

\section{The two dimensional case} \label{AdS2}

\subsection{Singularity and physical space}

The two dimensional case is very special because it doesn't
present a black hole structure. The particular structure directly
appears in the groupal formalism\footnote{See section \ref{BTZ} for notations related to $SL_2(\eR)$.}. Here
$G=SL_2(\eR)$ and, as homogeneous space, up to a double covering,

\begin{equation}
AdS_2=G/A=\Ad(G)H
\end{equation}
where $A=e^{\eR H}$ is the abelian part of $G$ with respect to the
Iwasawa decomposition. In the basis $\{H,E,F\}$ of
$SL_2(\eR)$, the matrix of the Killing form is given by

\begin{equation}
B=
\begin{pmatrix}
8&&\\
&&4\\
&4&
\end{pmatrix}
\end{equation}
while the basis  $\{H,E+F,E-F\}$ gives

\[
B=
\begin{pmatrix}
8\\
&8\\
&&-8
\end{pmatrix},
\]
so that we have the following isometry,
$(\sldr,B)\sim(\eR^3,\eta_{1,2})$. It will be convenient to see $AdS_2$
as an hyperboloid in $\eR^3$. We will use the Cartan involution $\theta(X)=-X^t$.

>From Definition \ref{Singular}, the singularities are here the closed orbits of $AN$
and $A\overline{N}$ for the adjoint action on $AdS_2=\Ad(G)H$. A
basis of the Lie algebra $\sA\oplus\sN$ is given by $\{E,H\}$. So
$x$ will belong to a closed orbit if and only if $E_x^*\wedge
H^*_x=0$. If we put $x=x_HH+x_EE+x_FF$, the computation is

\[
\begin{split}
E_x^*\wedge H^*_x&=[E,x]\wedge[H,x]\\
                 &=4x_Hx_F E\wedge F+2x_Ex_F H\wedge E-2x_F^2 H\wedge F.
\end{split}
\]
It is zero if and only if $x_F=0$. The closed orbit of
$A\overline{N}$ is given by the same computation with $H^*_x\wedge
F^*_x$. The part of these orbits contained in $AdS_2$ is the one
with norm $8$ :

\[
B(x,x)=8(x_H^2+x_Ex_F).
\]
In both cases, it gives $x_H=\pm 1$, and the closed orbits in
$AdS_2$ are given by

\begin{subequations}
\begin{align}
\pm H&+\lambda F\\
\pm H&+\lambda E,
\end{align}
\end{subequations}

\begin{proposition}\label{HSing}
 The singularities can equivalently be defined as
\begin{equation}\label{condHH}
\hS=\{x \in \Ad(G)H \tq \|H^*_x\|=0\}
\end{equation}
where $H^*$ is the fundamental field associated to the vector $H$:

\begin{equation}
H^*_x=\Dsdd{x\cdot e^{-tH}}{t}{0} =\Dsdd{\Ad(e^{-tH})x}{t}{0}
=-[H,x].
\end{equation}
\end{proposition}

\begin{proof}
The condition \eqref{condHH} for $x$ to belong to the singularity is

\begin{equation}\label{eq:BHxHx}
B([H,x],[H,x])=0.
\end{equation}
The most general\footnote{It is actually \emph{more} than the most
general element to be considered because our space is $\Ad(G)H$,
and not the whole $\sldr$.} element $x$ in $\sldr$ is
$x=aH+bE+cF$. It is easy to see that $[x,H]=-2bE+2cF$, so that the
condition \eqref{eq:BHxHx} becomes $bc=0$. Then the two
possibilities are $x=aH+bE$ and $x=aH+cF$. The singularities in
$\sldr$ are the planes $(H,F)$ and $(H,E)$. The intersection
between the plane $(H,F)$ and the hyperboloid is given by the
equation

\[
B(aH+bF,aH+bF)=8
\]
whose solutions are $a=\pm 1$. The same is also true for the plane
$(H,E)$. So we find back the fact that the singularities are given by the four lines
\begin{equation}
 \pm H+\lambda E \textrm{ and } \pm H+\lambda F.
\end{equation}
\end{proof}

Another way to express the singularities is
\begin{equation}
\Ad(e^{nE})(\pm H) \textrm{ and }
\Ad(e^{fF})(\pm H),
\end{equation}
which clearly shows that these are orbits of $AN$ and $A\overline{N}$. Indeed, as
$\Ad(a)$ fixes $H$, we can write $\Ad(an)H=\Ad(ana\emu)H$. Using
the CBH formula we find

\[
ana\emu=e^{nE+2anE+\ldots}=e^{ne^{2a}E}=n'\in N.
\]
The same can be done with $f$. So $\Ad(an)H=\Ad(n')H$ and
$\Ad(af)H=\Ad(f')H$. This shows that for all $n\in N$ and $a\in
A$, there exists a $n'\in N$ such that

\begin{subequations} \label{eq:singuAdd}
\begin{align}
\Ad(an)H&=\Ad(n')H \intertext{ The same is true with $f$ :}
\Ad(af)H&=\Ad(f')H.
\end{align}
\end{subequations}

In the basis $E,F,H$ the singularities are four lines with
angle=45\textdegree{} trough $H$ and $-H$. They divide the space
$AdS_2$ into four pieces. We define the \defe{physical space} as
the part of $AdS_2$ contained between $H+\lambda E$ and
$-H+\lambda E$.   The $K$ part of $SL_2(\eR)$ gives a double
covering of this curve. The part contained between the
singularities $H+\lambda F$ and $-H+\lambda F$ should be another
choice of physical space.

The following proposition gives an useful characterization of the
physical space.

\begin{proposition} \label{prop:AdAK}
Any point in the physical space can be written as $\Ad(ak)H$, with
$k\in]0,\pi/2[$.
\end{proposition}

\begin{proof}
The physical space contains the curve  $\cos\beta
H+\sin\beta(E+F)$ with $\beta\in]0,\pi[$, which is exactly
$\Ad(k)H$ for $k\in]0,\pi/2[$. It is also the intersection of
$AdS_2$ and the part of $\sldr$ between the planes $(E,H)$ and
$(F,H)$. If we use the coordinates $x,y,z$ on $\sldr$ (i.e.
$\overline{r}=xH+yE+zF$), our physical space is given by the
inequations

\[
\left\{
\begin{split}
x^2&+yz=1\\
y&>0\\
z&>0.
\end{split}
\right.
\]
The first equation gives a $\beta$ such that $x=\cos\beta$,
$yz=\sin^2\beta$. It is always possible to define a $a\in\eR$ such
that $y=e^{2a}\sin\beta$ and $z=e^{-2a}\sin\beta$. Finally, the
physical space is parameterized by

\begin{equation}
\overline{r}=\cos\beta H+\sin\beta(e^{2a}E+e^{-2a}F).
\end{equation}

On the other hand, from commutation relations in $\sldr$, one
finds

\begin{subequations}
\begin{align}
\Ad(e^{aH})E&=e^{2a}E,\\
\Ad(e^{aH})F&=e^{-2a}F.\\
\end{align}
\end{subequations}
Then
\begin{equation}
\begin{split}
\Ad(ak)H&=\Ad(e^{aH})(\cos\beta H+\sin\beta(E+F) )\\
    &=\cos\beta H+\sin\beta(e^{2a}E+e^{-2a}F).
\end{split}
\end{equation}
\end{proof}

\subsection{Light cone}

The light-like vectors of $\sldr$ are $E$ and $F$, so at
$\Ad(g)H$, the light-cone consists in two parts :

\[
\Ad(g)\Ad(e^{tE})H\text{ and } \Ad(g)\Ad(e^{tF})H .
\]
It is best rewritten in the compact form

\begin{equation}
C^+_{\Ad(g)H}=\{\Ad(g)\Ad(e^{t\epsilon E})H\}_{%
\begin{subarray}{l}
t>0\\
\epsilon=\id,\theta
\end{subarray}
}
\end{equation}
where $\epsilon$ is the identity or the Cartan involution.

It is somewhat easy to remark that for all $X,Y$ in a Lie algebra and all automorphism $\varphi$, the formula
$\varphi(\Ad(e^X)Y)=\Ad(e^{\varphi X})(\varphi Y)$ holds. Then

\begin{equation}
\Ad(e^{t\epsilon E})H=s(\epsilon)\epsilon(\Ad(e^{tE})H)
\end{equation}
with

\[
  s(\epsilon)=
\begin{cases}
1&\text{if $\epsilon=\id$},\\
-1&\text{if $\epsilon=\theta$}.
\end{cases}
\]
Since $H^*_x=-[H,x]$, the intersection of the light-cone with the singularity is expressed,
 using Proposition \ref{HSing}, as

\begin{equation}
\|[H,\Ad(g)\Ad(e^{t\epsilon E})H]\|^2=0.
\end{equation}

\subsection{No black hole}

The light cone of the point $\Ad(ak)H$ --which is a general point
of the physical space-- is given by
$\Ad(ak)s(\epsilon)\epsilon(\Ad(e^{tE})H)$. The computation of
$\Ad(ak)(H-2tE)$ and $-\Ad(ak)(-H+2tF)$ gives

\begin{subequations}  \label{eq:lin_light}
\begin{align}
(\cos(2k)-t\sin(2k))H-e^{2a}(\sin(2k)+2t\cos^{2}k)E-e^{-2a}(\sin(2k)-2t\sin^2k)F
\label{eq:lin_light_a} \intertext{and}
(\cos(2k)-t\sin(2k))H-e^{2a}(\sin(2k)-2t\sin^{2}k)E-e^{-2a}(\sin(2k)+2t\cos^2k)F
\end{align}
\end{subequations}
With respect to $t$, these are two straight lines, so they are the intersection of
$AdS_2$ and the tangent plane to $AdS_2$ at $\Ad(ak)H$.

This is important because it allows us immediately to infer
the non-existence of a black-hole structure for this choice of
singularity. The light cone at $x\in AdS_2$ is given by the
tangent plane $C$ of $AdS_2$ at $x$. The part of the singularity
passing by $H$ is given by a vertical plane $S$. The intersection
of these two planes is a line, and the intersection of a line with
$AdS_2$ is two points. Then each of the two lines of $C\cap AdS_2$
intersect one of the two lines of $S\cap AdS_2$. The same is true
for the other part of the singularity.

The conclusion is that both two lines of the light cone intersect
the singularity passing by $H$ \emph{and} the one passing by $-H$.
So any point comes from the singularity and returns to the
singularity; no point is connected to the infinity.

\section{Explicit matrix choices}\label{app_calc}

The first choice is to parameterize $SO(2,n)$ and $SO(1,n)$ in
such a way the latter leaves unchanged the vector
$(1,0,0,\ldots)$. Then

\begin{equation}  \label{eq:mtrH}
\sH=\soun\leadsto
  \begin{pmatrix}
     \begin{matrix}
       0&0\\
       0&0
     \end{matrix}
                       &  \begin{pmatrix}
                     \cdots 0\cdots\\
                \leftarrow v^t\rightarrow
                          \end{pmatrix}\\
    \begin{pmatrix}
       \vdots & \uparrow\\
         0    & v \\
       \vdots & \downarrow
    \end{pmatrix} &  B
  \end{pmatrix}.
\end{equation}
where  $v$ is $n\times 1$ and $B$ is skew symmetric $n\times n$.
 When we speak about $\mathfrak{so}(n)$, we usually refer to $B$.
 A complementary space $\sQ$ such that $[\sH,\sQ]\subset\sQ$ is given by

\begin{equation}
\sQ\leadsto
 \begin{pmatrix}
     \begin{matrix}
       0&a\\
       -a&0
     \end{matrix}
                       &  \begin{pmatrix}
              \leftarrow w^t\rightarrow \\
                 \cdots 0\cdots\\
                          \end{pmatrix}\\
    \begin{pmatrix}
      \uparrow   & \vdots\\
          w      &  0\\
      \downarrow & \vdots
    \end{pmatrix} & 0
  \end{pmatrix}
\end{equation}
We consider the involutive automorphism $\sigma=\id_{\sH}\oplus(-\id)_{\sQ}$ and the corresponding symmetric space structure on $\sG$. As basis of $\sQ$, we choice $q_0$ as the $2\times 2$ antisymmetric upper-left square and as $q_i$, the one obtained with $w$ full of zero apart a $1$ on the $i$th component.
Next we choice the Cartan involution $\theta(X)=-X^t$ which gives rise to a Cartan
 decomposition

\[
\sG=\sK\oplus\sP.
\]
The latter choice is made in such a way that $[\sigma,\theta]=0$.
It can be computed, but it is not astonishing that the compact part $\sK$ is made
 of  ``true'' rotations while $\sP$ contains the boost. So

\[
  \sK=
\begin{pmatrix}
  \sod\\
&\son
\end{pmatrix}.
\]
In order to build an Iwasawa decomposition, one has to choose a
maximal abelian
 subalgebra $\sA$ of $\sP$. Since rotations are in $\sK$, they must
  be boosts and the fact that there are only two time-like directions
  restricts
   $\sA$ to a two dimensional algebra. Up to reparametrization,
   it is thus generated by
 $t\partial_x+x\partial_t$ and $u\partial_y+y\partial_t$. Our matrix choices are

\[
   J_1=
\begin{pmatrix}
&0\\
0&0&0&1\\
&0\\
&1
\end{pmatrix}\in\sH,
\textrm{ and }
J_2=q_1=
\begin{pmatrix}
0&0&1&0\\
0\\
1\\
0
\end{pmatrix}\in\sQ.
\]
From here, we have to build root spaces. There still remains a lot
of arbitrary choices --among them, the positivity notion on the
dual space $\sA^*$. An elements $X$ in $\sG_{(a,b)}$ fulfill
$\ad(X)J_1=aJ_1$ and $\ad(X)J_2=bJ_2$. The symbol $E_{ij}$ denote
the matrix full of zeros with a $1$ on the component $ij$. Results
are

\begin{equation}
\sG_{(0,0)}\leadsto
\begin{pmatrix}
&&x&0\\
&&0&y\\
x&0\\
0&y\\
&&&& D
\end{pmatrix},
\end{equation}
where $D\in M_{(n-2)\times(n-2)}$ is skew-symmetric,

\begin{subequations}
\begin{align}
\sG_{(1,0)}&\leadsto W_i=E_{2i}+E_{4i}+E_{i2}-E_{i4},\\
\sG_{(-1,0)}&\leadsto Y_i=-E_{2i}+E_{4i}-E_{i2}-E_{i4},\\
\sG_{(0,1)}&\leadsto V_i=E_{1i}+E_{3i}+E_{i1}-E_{i3},\\
\sG_{(0,-1)}&\leadsto X_i=-E_{1i}+E_{3i}-E_{i1}-E_{i3}
\end{align}
\end{subequations}
with $\dpt{i}{5}{n+2}$ and

\begin{equation}
\sG_{(1,1)}\leadsto M=
\begin{pmatrix}
   0&1&0&-1\\
   -1&0&1&0\\
   0&1&0&-1\\
   -1&0&1&0
\end{pmatrix},
\quad
\sG_{(1,-1)}\leadsto L=
\begin{pmatrix}
   0&1&0&-1\\
   -1&0&-1&0\\
   0&-1&0&1\\
   -1&0&-1&0
\end{pmatrix},
\end{equation}

\begin{equation}
\sG_{(-1,1)}\leadsto N=
\begin{pmatrix}
   0&1&0&1\\
   -1&0&1&0\\
   0&1&0&1\\
   1&0&-1&0
\end{pmatrix},
\quad
\sG_{(-1,-1)}\leadsto F=
\begin{pmatrix}
   0&1&0&1\\
   -1&0&-1&0\\
   0&-1&0&-1\\
   1&0&1&0
\end{pmatrix}.
\end{equation}
The choice of positivity is

\begin{equation}
   \sN=\{V_i,W_j,M,L\}.
\end{equation}

The following result is important in the computation of the light
cones : if $k\in SO(n)$, then the choice $E=q_0+q_2$ of nilpotent
element in $\sQ$ gives

\begin{equation} \label{eq:AdkE}
   \Ad(k)E=
\begin{pmatrix}
0&1&w_1&w_2&\ldots\\
-1\\
w_1\\
w_2\\
\vdots
\end{pmatrix}
\end{equation}
where the vector $w$ is the first column of $k$, whose components
satisfy $\sum_{i=1}^{l-1} \, w_i^2=1$.

\addcontentsline{toc}{section}{References}


\begin{thebibliography}{10}

\bibitem{BTZ}
M.~Banados, C.~Teitelboim, and J.~Zanelli, {\it The black hole in
  three-dimensional space-time},  {\em Phys. Rev. Lett.} {\bf 69} (1992)
  1849--1851, [\href{http://xxx.lanl.gov/abs/hep-th/9204099}{{\tt
  hep-th/9204099}}].

\bibitem{BHTZ}
M.~Banados, M.~Henneaux, C.~Teitelboim, and J.~Zanelli, {\it
Geometry of the
  (2+1) black hole},  {\em Phys. Rev.} {\bf D48} (1993) 1506--1525,
  [\href{http://xxx.lanl.gov/abs/gr-qc/9302012}{{\tt gr-qc/9302012}}].

\bibitem{BDHRS}
P.~Bieliavsky, S.~Detournay, M.~Herquet, M.~Rooman, and
P.~Spindel, {\it Global
  geometry of the 2+1 rotating black hole},  {\em Phys. Lett.} {\bf B570}
  (2003) 231--240, [\href{http://xxx.lanl.gov/abs/hep-th/0306293}{{\tt
  hep-th/0306293}}].

\bibitem{Keio}
P.~Bieliavsky, S.~Detournay, P.~Spindel, and M.~Rooman, {\it
Noncommutative
  locally anti-de sitter black holes},
  \href{http://xxx.lanl.gov/abs/math.qa/0507157}{{\tt math.qa/0507157}}.

\bibitem{Figueroa}
J.~Figueroa-O'Farrill, O.~Madden, S.~F. Ross, and J.~Simon, {\it
Quotients of
  ads(p+1) x s**q: Causally well-behaved spaces and black holes},  {\em Phys.
  Rev.} {\bf D69} (2004) 124026,
  [\href{http://xxx.lanl.gov/abs/hep-th/0402094}{{\tt hep-th/0402094}}].

\bibitem{AdSBH}
M.~Banados, A.~Gomberoff, and C.~Martinez, {\it Anti-de sitter
space and black
  holes},  {\em Class. Quant. Grav.} {\bf 15} (1998) 3575--3598,
  [\href{http://xxx.lanl.gov/abs/hep-th/9805087}{{\tt hep-th/9805087}}].

\bibitem{Madden}
O.~Madden and S.~F. Ross, {\it Quotients of anti-de sitter space},
{\em Phys.
  Rev.} {\bf D70} (2004) 026002,
  [\href{http://xxx.lanl.gov/abs/hep-th/0401205}{{\tt hep-th/0401205}}].

\bibitem{Banados:1997df}
M.~Banados, {\it Constant curvature black holes},  {\em Phys.
Rev.} {\bf D57}
  (1998) 1068--1072, [\href{http://xxx.lanl.gov/abs/gr-qc/9703040}{{\tt
  gr-qc/9703040}}].

\bibitem{Aminneborg}
S.~Aminneborg, I.~Bengtsson, S.~Holst, and P.~Peldan, {\it Making
anti-de
  sitter black holes},  {\em Class. Quant. Grav.} {\bf 13} (1996) 2707--2714,
  [\href{http://xxx.lanl.gov/abs/gr-qc/9604005}{{\tt gr-qc/9604005}}].

\bibitem{HolstPeldan}
S.~Holst and P.~Peldan, {\it Black holes and causal structure in
anti-de sitter
  isometric spacetimes},  {\em Class. Quant. Grav.} {\bf 14} (1997) 3433--3452,
  [\href{http://xxx.lanl.gov/abs/gr-qc/9705067}{{\tt gr-qc/9705067}}].

\bibitem{BRS}
P.~Bieliavsky, M.~Rooman, and P.~Spindel, {\it Regular poisson
structures on
  massive non-rotating btz black holes},  {\em Nucl. Phys.} {\bf B645} (2002)
  349--364, [\href{http://xxx.lanl.gov/abs/hep-th/0206189}{{\tt
  hep-th/0206189}}].

\bibitem{BDRS}
P.~Bieliavsky, S.~Detournay, P.~Spindel, and M.~Rooman, {\it Star
products on
  extended massive non-rotating btz black holes},  {\em JHEP} {\bf 06} (2004)
  031, [\href{http://xxx.lanl.gov/abs/hep-th/0403257}{{\tt hep-th/0403257}}].

\bibitem{MemClem}
C.~Hyvrier, {\it Etude g\'eom\'etrique de certaines vari\'et\'es
localement
  anti-de sitter : les trous noirs {BTZ} massifs et sans moment angulaire},
  {\em M\'emoire de licence, Universit\'e Libre de Bruxelles}.

\bibitem{Helgason}
S.~Helgason, {\em Differential geometry and symmetric spaces}.
\newblock Pure and Applied Mathematics, Vol. XII. Academic Press, New York,
  1962.

\bibitem{Loos}
O.~Loos, {\em Symmetric spaces. {I}: {G}eneral theory}.
\newblock W. A. Benjamin, Inc., New York-Amsterdam, 1969.

\bibitem{kobayashi}
S.Kobayashi and K.Nomizu, {\em Foundation of differential
geometry}, vol.~1.
\newblock Interscience publishers, 1963.

\bibitem{kobayashi2}
S.Kobayashi and K.Nomizu, {\em Foundation of differential
geometry}, vol.~2.
\newblock Interscience publishers, 1969.

\bibitem{Kerin}
M.~Kerin and D.~Wraith, {\it Homogeneous metrics on spheres},
{\em Irish Math.
  Soc. Bull.} (2003), no.~51 59--71.

\end{thebibliography}
\providecommand{\href}[2]{#2}\begingroup\raggedright\endgroup

\end{document}